\documentclass{amsart}
\usepackage[english]{babel}
\usepackage[latin1]{inputenc}
\usepackage[dvips,final]{graphics}
\usepackage{amsmath,amsfonts,amssymb,amsthm,amscd,array,stmaryrd,mathrsfs}
\usepackage{pstricks}
 \usepackage[all]{xy}
\usepackage{textcomp}
 \usepackage[final]{epsfig}
\theoremstyle{plain}
\newtheorem{lem}{Lemma}[section]

\theoremstyle{definition}

\newtheorem{rem}[lem]{Remark}

\newcommand{\R}{\mathbb{R}}
\newcommand{\Z}{\mathbb{Z}}

\newcommand{\bbH}{\mathbb{H}}
\newcommand{\bbO}{\mathbb{O}}

\newcommand{\G}{\mathcal{G}}

\newcommand{\gog}{\mathfrak{g}}

\newcommand{\e}{\varepsilon}

\def\e{\varepsilon}

\def\G{\Gamma}

\def\s{\sigma}

\begin{document}

\title{Well, Papa, can you multiply triplets?}

\author{Sophie Morier-Genoud
\hskip 1cm
Valentin Ovsienko}

\address{
Sophie Morier-Genoud,
Universit\'e Paris Diderot Paris 7,
UFR de math\'ematiques case 7012,
75205 Paris Cedex 13, France}

\address{
Valentin Ovsienko,
CNRS,
Institut Camille Jordan,
Universit\'e Claude Bernard Lyon~1,
43 boulevard du 11 novembre 1918,
69622 Villeurbanne cedex,
France}

\email{sophiemg@umich.edu,
ovsienko@math.univ-lyon1.fr}

\date{}


\begin{abstract}
We show that the classical algebra of quaternions
is a commutative $\Z_2\times\Z_2\times\Z_2$-graded algebra.
A similar interpretation of the algebra of octonions is impossible.
\end{abstract}

\maketitle

\thispagestyle{empty}


This note is our ``private investigation'' of what really happened
on the 16th of October, 1843 on the Brougham Bridge
when Sir William Rowan Hamilton engraved
on a stone his fundamental relations:
$$
i^2=j^2=k^2=i\cdot{}j\cdot{}k=-1.
$$
Since then,
the elements $i,j$ and $k$, together with the unit, $1$,
denote the canonical basis of the celebrated
4-dimensional associative algebra of quaternions
$\bbH$.

Of course, the algebra $\bbH$ is not commutative: the relations above
imply that the elements $i,j,k$ \textit{anticommute} with each other, for instance
$$
i\cdot{}j=-j\cdot{}i=k.
$$
Yes, but...

\section{The algebra of quaternions is commutative}

Our starting point is the following amazing observation.

\bigskip

\textit{\centerline{The algebra $\bbH$ is, indeed, commutative.}}

\bigskip

Let us introduce the following ``triple degree'':
\begin{equation}
\label{DegH}
\begin{array}{rcl}
\s(1)&=&(0,0,0),\\[4pt]
\s(i)&=&(0,1,1),\\[4pt]
\s(j)&=&(1,0,1),\\[4pt]
\s(k)&=&(1,1,0),
\end{array}
\end{equation}
then, quite remarkably, the usual product of quaternions
satisfies the graded commutativity condition:
\begin{equation}
\label{ComPr}
p\cdot{}q=(-1)^{\left\langle\s(p),\s(q)\right\rangle}\,q\cdot{}p,
\end{equation}
where $p,q\in\bbH$ are \textit{homogeneous}
(i.e., proportional one of the basic vector) and where $\langle\,,\,\rangle$ 
is the scalar product of 3-vectors.
Indeed, $\left\langle\s(i),\s(j)\right\rangle=1$ and similarly for $k$,
so that $i,j$ and $k$ anticommute with each other, but
$\left\langle\s(i),\s(i)\right\rangle=2$.
The product $i\cdot{}i$ of $i$ with itself is commutative and similarly
for $j$ and $k$, without any contradiction.

The degree (\ref{DegH}) viewed as an element of
the abelian group $\Z_2\times\Z_2\times\Z_2$ satisfies
the following linearity condition
\begin{equation}
\label{Grading}
\s(x\cdot{}y)=\s(x)+\s(y),
\end{equation}
for all homogeneous $x,y\in\bbH$.
The relations (\ref{ComPr}) and (\ref{Grading}) together mean that
$\bbH$ is a $\Z_2\times\Z_2\times\Z_2$-graded commutative algebra.

We did not find the above observation in the literature
(see however \cite{Mag} for a different ``abelianization'' of
$\bbH$ in terms of a twisted $\Z_2\times\Z_2$ group algebra,
see also \cite{Bae,Con,Lam}).
Its main consequence is a systematic procedure of \textit{quaternionization}
(similar to complexification).
Indeed, many classes of algebras allow tensor product with commutative
algebras.
Let us give an example.
Given an arbitrary real Lie algebra $\gog$,
the tensor product $\gog_\bbH:=\bbH\otimes_\R\gog$ is a
 $\Z_2\times\Z_2\times\Z_2$-graded Lie algebra.
If furthermore $\gog$ is a real form of a simple complex Lie algebra, then
$\gog_\bbH$ is again simple.

The above observation gives a general idea to study
graded commutative algebras over the abelian group
$$
\G=\underbrace{\Z_2\times\cdots\times\Z_2}_{n\,\rm{times}}.
$$
One can show that, in some sense, this is the most general grading,
in the graded-commutative-algebra context, but we will not provide the details here.
Let us mention that graded commutative algebras are essentially studied
in the case $\G=\Z_2$ (or $\Z$),
almost nothing is known in the general case.

\section{...but not the algebra of octonions}

After the quaternions,
the next ``natural candidate'' for commutativity
is of course the algebra of octonions $\bbO$.
However, let us show that 

\bigskip

\textit{\centerline{The algebra $\bbO$
cannot be realized as a graded commutative algebra.}}

\bigskip

Indeed, recall that $\bbO$ contains 7 mutually anticommuting
elements $e_1,\ldots,e_7$ such that $(e_\ell)^2=-1$ for $\ell=1,\ldots,7$
that form several copies of $\bbH$
(see \cite{Bae,Con} for beautiful introduction to the octonions).
Assume there is a grading $\s:e_\ell\mapsto\G$
with values in an abelian group $\G$,
satisfying (\ref{ComPr}) and (\ref{Grading}).
Then, for three elements
$e_{\ell_1},e_{\ell_2},e_{\ell_3}\in\bbO$ such that
$e_{\ell_1}\cdot{}e_{\ell_2}=e_{\ell_3}$, one has
$$
\s(e_{\ell_3})=\s(e_{\ell_1})+\s(e_{\ell_2}).
$$
If now $e_{\ell_4}$ anticommutes with $e_{\ell_1}$ and $e_{\ell_2}$,
then $e_{\ell_4}$ has to commute with $e_{\ell_3}$ because
of the linearity of the scalar product.
This readily leads to a contradiction.

\section{Multiplying the triplets}

Let us now take another look at  the grading (\ref{DegH}).
It turns out that there is a simple way to reconstitute the whole
structure of $\bbH$ directly from this formula.

First of all, we rewrite the grading as follows:
\begin{equation}
\label{NewGrad}
\begin{array}{rcl}
1& \leftrightarrow&(0,\,0,\,0),\\[4pt]
i& \leftrightarrow&(0,\,1_2,1_3),\\[4pt]
j& \leftrightarrow&(1_1,0,\,1_3),\\[4pt]
k& \leftrightarrow&(1_1,1_2,0).
\end{array}
\end{equation}
Second of all, we define the rule for multiplication of triplets.
This multiplication is nothing but the usual operation in
$\Z_2\times\Z_2\times\Z_2$, i.e., the component-wise addition
(modulo 2), for instance,
$$
(1_1,0,0)\cdot(1_1,0,0)=(0,0,0),
\qquad
(1_1,0,0)\cdot(0,1_2,0)=(1_1,1_2,0),
$$
but with an important additional \textit{sign rule}.
Whenever we have to exchange ``left-to-right''
two units, $1_n$ and $1_m$ with $n>m$,
we put the ``$-$'' sign, for instance
$$
(0,1_2,0)\cdot(1_1,0,0)=-(1_1,1_2,0),
$$
since we exchanged $1_2$ and $1_1$.

One then has for the triplets in (\ref{NewGrad}):
$$
i\cdot{}j \leftrightarrow
(0,1_2,1_3)\cdot(1_1,0,1_3)=
(1_1,1_2,0) \leftrightarrow k,
$$
since the total number of exchanges is \textit{even} 
($1_2$ and $1_3$ were exchanged with $1_1$) and
$$
j\cdot{}i \leftrightarrow
(1_1,0,1_3) \cdot(0,1_2,1_3)=
-(1_1,1_2,0)
\leftrightarrow-k,
$$
since  the total number of exchanges is \textit{odd} ($1_3$ was exchanged with $1_2$).
In this way, one immediately recovers the complete multiplication table of $\bbH$.

\begin{rem}
The above realization is of course related to the embedding
of $\bbH$ into the associative algebra with 3 generators
$\e_1,\e_2,\e_3$ subject to the relations
$$
\e_n^2=1,
\qquad
\e_n\e_m=-\e_m\e_n,
\quad
\hbox{for}
\quad 
n\not=m.
$$
This embedding is given by 
$$
i\mapsto\e_2\e_3,
\qquad
j\mapsto\e_1\e_3,
\qquad
k\mapsto\e_1\e_2
$$
and is well-known.
\end{rem}

Everybody knows the famous story of Hamilton and his son
asking his father the same question every morning: 
``Well, Papa, can you multiply triplets?''
and always getting the same answer:
``No, I can only add and subtract them'',
with a sad shake of the head.
This story has now a happy end.
As we have just seen,
Hamilton did nothing but multiplied the triplets.
Or should we rather say added and subtracted them?

\bigskip

\noindent \textbf{Acknowledgments}.
We are grateful to S. Tabachnikov for helpful suggestions.


\end{document}